\documentclass[12pt]{amsart}

\usepackage[margin=1in,letterpaper,portrait]{geometry}
\usepackage{amsmath}
\usepackage{amssymb}
\usepackage{amsthm}
\usepackage{mathrsfs}
\usepackage{verbatim}
\usepackage{fancyhdr}
\usepackage[pdftex,bookmarksnumbered=true,linkbordercolor={1 0 0}]{hyperref}
\usepackage{cleveref}
\usepackage{url}
\usepackage{tikz}
\usepackage{enumerate}
\usepackage{longtable}
\usetikzlibrary{matrix,arrows}
\usepackage{tikz-cd}
\usetikzlibrary{plotmarks}
\usepackage{setspace}
\newcommand{\student}{Dennis Tseng}
\newcommand{\studentemail}{DennisCTseng@gmail.com}

\newcommand{\mb}{\mathbb}

\newcommand{\mc}{\mathcal}
\newcommand{\im}{{\rm im}}
\newcommand{\ol}{\overline}
\newcommand{\wt}{\widetilde}
\newcommand{\codim}{{\rm codim}}
\newcommand{\Hilb}{{\rm Hilb}}
\newcommand{\RHilb}{\wt{{\rm Hilb}}}

\renewcommand{\O}{\mathscr{O}}

\newcommand{\Spec}{{\rm Spec}}

\newcommand{\rank}{{\rm rank}}

\newcommand*{\sheafhom}{\mathscr{H}\kern -.5pt om}
\newcommand{\on}[1]{\operatorname{#1}}

\newtheorem{theorem}{Theorem}[section]
\crefname{theorem}{Theorem}{Theorems}
\newtheorem{corollary}[theorem]{Corollary}

\newtheorem{problem}[theorem]{Problem}

\newtheorem{lemma}[theorem]{Lemma}
\newtheorem{claim}[theorem]{Claim}
\newtheorem{fact}[theorem]{Fact}
\newtheorem{proposition}[theorem]{Proposition}
\crefname{proposition}{Proposition}{Propositions}

\theoremstyle{definition}
\newtheorem{definition}{Definition}[section]

\newtheorem{remark}{Remark}

\title{On degenerate sections of vector bundles}
\author{Dennis Tseng}
\address{Dennis Tseng, Harvard University, Cambridge, MA 02138}
\email{DennisCTseng@gmail.com}
\date{\today}

\begin{document}
\footnotetext[1]{This material is based upon work supported by the National Science Foundation Graduate Research Fellowship Program under Grant No. 1144152.}
\maketitle
\begin{abstract}
We consider the locus of sections of a vector bundle on a projective scheme that vanish in higher dimension than expected. We will find the largest components of this locus asymptotically, after applying a high enough twist to the vector bundle. We will also give an interpretation in terms of a limit in the Grothendieck ring of varieties. 
\end{abstract}
\section{Introduction}
We will work throughout over an algebraically closed field $K$ of arbitrary characteristic. Let $X$ be an irreducible projective scheme and $V$ be a globally generated vector bundle on $X$. A general section $s\in H^0(V)$ either vanishes in codimension exactly $\rank(V)$ or is nonvanishing. Considering $H^0(V)$ as an affine space, there is a closed locus $D(V)\subset H^0(V)$ consisting of sections that vanish in codimension less than $\rank(V)$. We are interested in basic questions about this locus, for example:

\begin{problem}
\label{Q1}
What is the dimension of $D(V)$? What are the components, and what can we say about them? 
\end{problem}

The purpose of this paper is to give a clean answer after a twist of $V$ by a high tensor power of $\O_X(1)$. For example, we show
\begin{theorem}
\label{DIM}
There exists $N_0$ such that for all $N\geq N_0$, any section $s\in H^0(V(N))$ in any component of $D(V(N))$ of largest dimension vanishes on some subscheme of $X$ of codimension $\rank(V)-1$ and of minimal degree among all subschemes of codimension $\rank(V)-1$ in $X$.
\end{theorem}

Given \Cref{DIM}, one may then ask for next largest components of $D(V(N))$ and expect that they are given by sections that vanish on some subscheme of second largest degree among all varieties of codimension $\rank(V)-1$. Then, the next largest should correspond to varieties of third largest degree and so on. 

We will show that this is in fact true (see \Cref{MT}) and will use the Hilbert polynomial instead of degree as a finer invariant to distinguish between the dimensions of the various components of $D(V(N))$. In particular, the number of components of $D(V(N))$ will necessarily approach infinity as $N\to\infty$, but we will only have fine control of the dimensions of its large dimensional components. 

Finally, one can define $D(V(N),a)$ as the sections that vanish in codimension $\rank(V)-a$ for integer $a\geq 1$, so $D(V(N),1)=D(V(N))$. Even though the case $a=1$ is the most natural, our results hold for all $a\geq 1$. 

\subsection{Special case: excess intersection of hypersurfaces in projective space}
If $X=\mb{P}^r$ and $V=\O(d_1)\oplus \cdots \oplus \O(d_k)$, then \Cref{DIM} specializes to \Cref{HYP}, a statement about hypersurfaces failing to intersect properly in projective space.
\begin{corollary}
\label{HYP}
Let $d_1,\ldots,d_k$ be degrees. Then, there exists $N_0$ such that for all $N\geq N_0$, the locus $D(\O_{\mb{P}^r}(d_1+N)\oplus\cdots\oplus \O_{\mb{P}^r}(d_k+N))$ has a unique component of maximal dimension. For every tuple $(F_1,\ldots,F_k)$ of homogenous forms in that component, $\{F_1=\cdots=F_k=0\}$ contains some linear space of dimension $r-k+1$. 
\end{corollary}
Uniqueness is part of \Cref{MT}, and it also follows from an incidence correspondence, expressing the maximal component as the image of
\begin{align*}
    \prod_{i=1}^{k}H^0(\O_{\mb{P}^r}(d_i))\times \mb{G}(r-k+1,r)\supset \{(F_1,\ldots,F_k),\Lambda\mid F_1|_\Lambda=\cdots=F_k|_\Lambda=0\}\to  \prod_{i=1}^{k}H^0(\O_{\mb{P}^r}(d_i))
\end{align*}
given by forgetting the plane $\Lambda$. 

There are simple examples where the conclusion of \Cref{HYP} fails when $N=0$. For example, consider the case $X=\mb{P}^4$ and $V=\O_{\mb{P}^4}(2)\oplus \O_{\mb{P}^4}(2)$. Then, 
\begin{align*}
    \{(F_1,F_2)\mid \{F_1=F_2=0\}\text{ contains some hyperplane}\}
\end{align*}
is codimension 16 in $V$, as it is codimension 6 for each quadric to be reducible and a further codimension 4 for the reducible quadrics to share a hyperplane. On the other hand, 
\begin{align*}
    \{(F_1,F_2)\mid \text{$F_1$ and $F_2$ are related by a scalar multiple}\}
\end{align*}
is codimension 14 in $V$. 
In an effort to remove the large twist, the author previously obtained a quantitative version of Corollary \ref{HYP}, where the conclusion of \Cref{HYP} holds when $N=0$ if $k=r$ and the degrees $d_1,\ldots,d_r$ are not too far apart \cite[Theorem 1.3]{Tseng}. See also \cite{TRockyCIH} for a weaker statement that holds for all degrees $d_1,\ldots,d_r$ when $N=0$ and $k=r$. 

\subsection{Acknowledgements}
The author would like to thank Anand Patel for suggesting the problem and for helpful conversations on initial approaches. The author would also like to thank the anonymous referee during the revision of this paper. 

\section{Definitions and setup}
\label{defsec}
Throughout the paper, we will fix an irreducible projective scheme $X$ and a vector bundle $V$ on $X$. We will work over an algebraically closed field $K$ of arbitrary characteristic. Unless otherwise specified, $a$ is a positive integer. Given a scheme $Z$, we will denote by $Z^{\operatorname{red}}$ the reduction of $Z$ \cite[Proposition 4.2(c)]{Liu}.

There is a table of symbols and definitions in \Cref{symboltable} at the end of the document to help the reader.

\subsection{Hilbert schemes}
\begin{definition}
\label{hilbdef}
Let $\Hilb_X$ be the Hilbert scheme of subschemes of $X$ and $\RHilb_X\subset \Hilb_X$ denote the locus parameterizing geometrically irreducible subschemes, which is open \cite[IV 12.2.1(x)]{EGA}. To indicate dimension, we let $\Hilb_X^c$ and $\RHilb_X^c$ denote the restriction to the connected components parameterizing subschemes of dimension $c$. 
\end{definition}

Recall the Euler characteristic of a twist $F(t)$ of a coherent sheaf $F$ on a projective scheme 
\begin{align*}
    \chi(F(t)):=\sum_{i\geq 0}{(-1)^ih^i(F(t))}
\end{align*} 
is a polynomial in $t$, which is the Hilbert polynomial of $F$. From the usual proof of the polynomiality of the Hilbert polynomial in the case $F$ is the structure sheaf of a subscheme via considering hyperplane slices and applying additivity of Euler characteristic in short exact sequences, one has the following fact.

\begin{fact}
\label{VS}
Suppose $Z$ is an integral projective scheme and $F$ is a coherent sheaf on $Z$ of positive rank. Then, the Hilbert polynomial $\chi(F(t))$ is a polynomial of degree $\dim(Z)$ in $t$ with leading coefficient $\frac{\deg(Z)\rank(F)}{\dim(Z)!}$. 
\end{fact}

\begin{definition}
\label{hilbdefP}
Given a polynomial $P(t)\in \mb{Q}[t]$, let $\RHilb_X^{c}(V,P(t))\subset \RHilb_X^{c}$ denote the connected components of $\RHilb_X^{c}$ parameterizing subschemes $Z$ where $\chi(V|_{Z}(t))=P(t)$. 
\end{definition}

\begin{definition}
\label{hilbdefPc}
Let $\ol{\RHilb_X^{c}(V,P(t))}$ be the closure of $\RHilb_X^{c}(V,P(t))$ in $\Hilb_X$. 
\end{definition}

\subsection{Hilbert polynomials and Hilbert functions}

\begin{definition}
\label{Sc}
Let $S_c$ be the set of polynomials $P(t)\in \mb{Q}[t]$ such that $\RHilb_X^{c}(V,P(t))$ is nonempty.
\end{definition}
Necessarily, $S_c$ consists of polynomials of degree $c$ (see \Cref{VS}). 
\begin{definition}
\label{pV}
Let $P_V(t)\in \mb{Q}[t]$ denote the Hilbert polynomial $\chi(X,V(t))$. 
\end{definition}

\begin{definition}
\label{HIF}
Given $Z\subset X$ and a vector bundle $E$ on $X$, define the Hilbert function $h_{Z,E}$ to be
\begin{align*}
h_{Z,E}(n):= \dim(\im(H^0(E(n))\rightarrow H^0(E(n)|_{Z}))).
\end{align*}
\end{definition}
In particular, Definition \ref{HIF} depends on our ambient projective scheme $X$. 

\begin{definition}
\label{PQ}
Given two polynomials $P(t),Q(t)\in \mb{Q}[t]$, we say $P(t)$ \emph{dominates} $Q(t)$ if $\lim_{t\rightarrow\infty}P(t)-Q(t)=\infty$ and $P(t)$ is \emph{equivalent} to $Q(t)$ if neither $P(t)$ nor $Q(t)$ dominates the other (e.g. $P(t)$ and $Q(t)$ differ by a constant). If $P(t)$ and $Q(t)$ are equivalent, we will also write this as $P(t)\sim Q(t)$.
\end{definition}

Recall Chow's finiteness theorem.
\begin{theorem}[{\cite[Exercise I.3.28 and Theorem I.6.3]{Kollar}}]
\label{Chow}
There are only finitely many components of $\RHilb_X$ parameterizing subschemes of a fixed dimension and degree.
\end{theorem} 
This means there exists a sequence $P_{c,1},P_{c,2},\ldots \in \mb{Q}[t]$ such that $\RHilb_X^{c}(V,P_{c,i}(t))$ is nonempty for each $i$, $P_{c,i+1}$ dominates $P_{c,i}$ for each $i$, and for every $P$ for which $\RHilb_X^{c}(V,P(t))$ is nonempty, $P(t)\sim P_{c,i}(t)$ for some $i$. 

\begin{definition}
\label{Pi}
Let $P_{c,1}(t),P_{c,2}(t),\ldots$ be a fixed choice of polynomials in $\mb{Q}[t]$ with the properties above. 
\end{definition}

\subsection{Degenerate sections of a vector bundle}
\begin{definition}
\label{DA}
Let $D(V,a)\subset H^0(V)$ denote the closed locus of sections $s\in H^0(V)$ such that $\{s=0\}$ has codimension at most $\rank(V)-a$ in $X$. 
\end{definition}

\begin{definition}
\label{DI}
Let $\wt{D}(V, a, P(t))$ consist of pairs $(s,[Z])\in H^0(V)\times \RHilb_X^{\dim(X)-\rank(V)+a}(V,P(t))$, where $s$ vanishes on $Z$. The locus $\wt{D}(V,a,P(t))\subset H^0(V)\times \RHilb_X^{\dim(X)-\rank(V)+a}(V,p(t))$ is closed.
%and it can even be given a canonical scheme structure as an open subset of a relative Hilbert scheme \cite[Lemma 7.1]{ACGH2}.
\end{definition}

We have the incidence correspondence
\begin{center}
\begin{tikzcd}
& \wt{D}(V,a,P(t)) \arrow[dl,  "\pi_1" , swap] \arrow[dr, "\pi_2"] &\\
 D(V,a) & &  \RHilb_X^{\dim(X)-\rank(V)+a}(V,a,P(t))
\end{tikzcd}
\end{center}

\begin{definition}
\label{DIp}
Given $P(t)\in \mb{Q}[t]$, let $D(V,a,P(t))$ be the constructible subset $\pi_1(\wt{D}(V,a,P(t)))$ of $H^0(V)$. 
\end{definition}

\begin{definition}
\label{DIc}
Like in Definition \ref{DI}, let $\wt{D}(V, a, P(t))^{\text{cl}}$ denote the pairs $(s,[Z])\in H^0(V)\times \ol{\RHilb_X^{\dim(X)-\rank(V)+a}(V,P(t))}$  where $s$ vanishes on $Z$. 
\begin{center}
\begin{tikzcd}
& \wt{D}(V,a,P(t))^{\text{cl}} \arrow[dl,  "\pi_1" , swap] \arrow[dr, "\pi_2"] &\\
 D(V,a) & &  \ol{\RHilb_X^{\dim(X)-\rank(V)+a}(V,a,P(t))}
\end{tikzcd}
\end{center}
Let $D(V,a,P(t))^{\text{cl}}$ be the closed subset $\pi_1(\wt{D}(V, a, P(t))^{\text{cl}})$ of $H^0(V)$. 
\end{definition}

\subsection{Constructible sets}
\label{CONSET}
Since we will need to work with constructible sets, such as $D(V,a,P(t))$ in Definition \ref{DIp}, we will recall basic facts about dimensions of constructible sets. To take the dimension of a constructible set, it suffices to either look at the generic points or take the closure.
\begin{definition}
\label{DEFDIM}
If $A\subset X$ is a constructible set, then $\dim(A):=\dim(\overline{A})$.
\end{definition}

\begin{lemma}
\label{FY6}
If $f:X\rightarrow Y$ is a morphism of finite type $K$-schemes, $A\subset X$ and $f(A)=B\subset Y$ constructible sets, and $\dim(f^{-1}(b)\cap A)<c$ for all $b\in B$, then
\begin{align*}
\dim(A)\leq \dim(B)+c.
\end{align*}
If $\dim(f^{-1}(b))=c$ for all $b\in B$, then equality holds.
\end{lemma}
\begin{proof}
Apply the usual theorem on fiber dimension at the generic points of the components of $\overline{A}$ to $f|_{\ol{A}}: \overline{A}\rightarrow \overline{B}$. 
\end{proof}

\begin{definition}
If $A\subset B$ are constructible subsets of a scheme $Y$, then the codimension of $A$ in $B$ is defined to be $\dim(B)-\dim(A)$. If $A$ is empty, then the codimension is $\infty$. 
\end{definition}

\section{Statement of Results}
\label{statementsec}
Given the notation in \Cref{defsec}, we can now state our results. First, \Cref{MT} says the largest components of $D(V(N),a)$ consist of sections that vanish on subschemes of dimension $\dim(X)-\rank(V)+a$ of small Hilbert polynomial. Throughout this section, let $c=\dim(X)-\rank(V)+a$ and $a>0$ a positive integer. 

\begin{theorem}
\label{MT}
For each $m\geq 0$ and $a>0$, there exists $N_0$ such that for $N\geq N_0$, (the closure of) a component of maximal dimension of 
\begin{align*}
D(V(N),a)\backslash \bigcup_{i=1}^{m}\bigcup_{P(t)\sim P_{c,i}(t)}D(V(N),a,P(t))
\end{align*}
is a component of $\ol{D(V(N),a,P(t))}$ for some $P(t)\sim P_{c,m+1}(t)$.
\end{theorem}

In particular, setting $m=0$ identifies the largest components of $D(V(N),a)$ as sections vanishing on subschemes of dimension $c$ with minimal Hilbert polynomial. From the proof, we will also see that the components of $D(V(N),a)$ of largest dimensions look like vector bundles over components of $\RHilb_X^{c}$ away from subsets of large codimension. This is stated in \Cref{incidence} below. 

\begin{theorem}
\label{incidence}
For each $m\geq 0$ and $a>0$, there exists $N_0$ such that for $N\geq N_0$, we have an open subset $U_N$ of
\begin{align}
    D(V(N),a)\backslash \bigcup_{i=1}^{m}\bigcup_{P(t)\sim P_{c,i}(t)}D(V(N),a,P(t)) \label{subtract}
\end{align}
that is also an open subset of the vector bundle 
\begin{align*}
    \mathcal{V}=\left(\bigcup_{P(t)\sim P_{c,m+1}(t)}\wt{D}(V, a, P(t))\right)^{\on{red}}\to \left(\bigcup_{P(t)\sim P_{c,m+1}(t)}\RHilb_{X}^{\dim(X)-\rank(V)+a}(V,P(t))\right)^{\on{red}},
\end{align*} 
which is of rank $P_V(N)-P(N)$ on each component $\RHilb_{X}^{\dim(X)-\rank(V)+a}(V,P(t))$. Furthermore, the codimension of the complement of $U_N$ in \eqref{subtract} and in $\mathcal{V}$ tends to $\infty$ as $N\to \infty$.
\end{theorem}
\begin{center}
    \includegraphics[scale=.4]{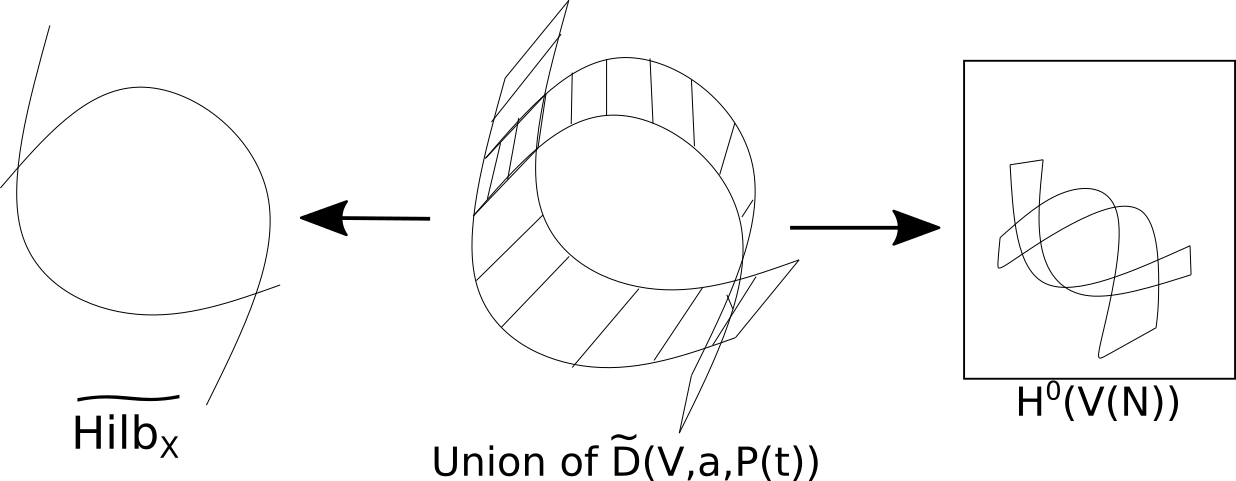}
\end{center}
\subsection{Interpretation with the Grothendieck Ring}
\label{GRdef}
We interpret \Cref{incidence} using the Grothendieck ring of varieties in \Cref{GR} below. This subsection is separate from the rest of the paper.

As an abelian group, the Grothendieck ring $\mc{M}$ is generated by the classes of finite type $K$-schemes up to isomorphism. Given a $K$-scheme $Z$, the class in the Grothendieck ring is denoted by $[Z]$ and $\bf{L}$ is defined to be $[\mb{A}^1]$. The group relations in the Grothendieck ring are generated by the scissor relations
\begin{align*}
    [Z]=[Z\backslash W]+[W]
\end{align*}
for $W\subset Z$ a closed subscheme. In particular, letting $W=Z^{\text{red}}$ shows $[Z]=[Z^{\text{red}}]$, so the Grothendieck ring does not see scheme-theoretic structure. Products in $\mc{M}$ are given by products over $K$ so $[Z_1]\cdot [Z_2]=[Z_1\times Z_2]$. The unit in $\mc{M}$ is $[\operatorname{Spec}(K)]$. 

See \cite{CNS11} for a survey on Grothendieck rings. However, we will mainly rely on the following fact, which follows from the scissor relations and Noetherian induction.
\begin{fact}
\label{affinebundle}
Let $Z$ be a finite type $K$-scheme and $W\to Z$ be a Zariski-trivial $\mathbb{A}^r$-bundle. Then, $[W]=[Z]\bf{L}^r$ in $\mc{M}$.
\end{fact}

Let $\mc{M}_{\bf{L}}$ be $\mc{M}$ with $\bf{L}$ inverted. We have a filtration 
\begin{align*}
\cdots\subset F_{-1}\subset F_{0}\subset F_{1}\subset \cdots \subset \mc{M}_{\bf{L}}
\end{align*}
given by dimension, and we can complete to obtain a ring $\widehat{\mc{M}_{\bf{L}}}$ \cite[Section 1.5]{Bourqui}. For more details on the Grothendieck ring and another example of a limit in $\widehat{\mc{M}_{\bf{L}}}$, see \cite{VakilWood}.

\Cref{GR} interprets \Cref{incidence} in terms of the Grothendieck ring and follows directly from its statement.
\begin{theorem}
\label{GR}
In $\widehat{\mc{M}_{\bf{L}}}$, for $m\geq 0$ and $a>0$,
\begin{align*}
\lim_{N\to \infty}\frac{[D(V(N),a)\backslash \bigcup_{i=1}^{m}\bigcup_{p(t)\sim p_i(t)}D(V(N),a,P(t))]}{{\bf L}^{P_V(N)-P_{c,m+1}(N)}} = \sum_{P(t)\sim P_{c,m+1}(t)}\frac{[\RHilb_{X}^{\dim(X)-\rank(V)+a}(V,P(t))]}{{\bf L}^{P(N)-P_{c,m+1}(N)}}. 
\end{align*}
\end{theorem}
The set of all $P(t)$ satisfying $P(t)\sim P_{c,m+1}(t)$ is finite by Chow's finiteness theorem (\Cref{Chow}) and \Cref{VS}.

\section{Proofs}
\label{proofs}
In this section, we prove \Cref{MT,incidence}. In the proofs, we can and will assume that $V$ is 0-regular in the sense of Castelnuovo-Mumford regularity, namely that $H^i(V(-i))=0$ for all $i> 0$ \cite[Definition 1.8.1]{PositivityI}. We will mainly use the consequence that $V$ is globally generated \cite[Theorem 1.8.3]{PositivityI}. We will use the full assumption of 0-regularity rather than just global generation exactly once in the proof of \Cref{TT1}. As in \Cref{statementsec}, we set $c=\dim(X)-\rank(V)+a$ and assume $a\geq 1$.
\subsection{Case of low degree}
\label{LOD}
In \Cref{LOD}, we will consider the components $\wt{D}(V(N),a,P(t))$, where the degree of $P(t)$ is bounded and $N$ is increasing. The strategy is to set up the usual incidence correspondence. 
\begin{proposition}
\label{VBN}
Given $P(t)\in S_c$, there exists $N_0$ dependent on $P(t)$ such that for all $N\geq N_0$, 
\begin{align*}
\left(\wt{D}(V(N),a,P(t))^{\text{cl}}\right)^{\operatorname{red}}\rightarrow \left(\ol{\RHilb_X^{c}(V,P(t))}\right)^{\operatorname{red}}
\end{align*}
is a vector bundle of rank $P_V(N)-P(N)$. 
\end{proposition}
We must have $P_V(N)\geq P(N)$ for $N$ large enough because $V|_{Z}$ is a quotient of $V$ for subschemes $Z$ of $X$. 

\begin{proof}[Proof of \Cref{VBN}]
\begin{claim}
We can choose $N_0$ large enough so that for each $[Z]\in \ol{\RHilb_X^{c}(V,P(t))}$, the Hilbert function and Hilbert polynomial agree for $N\geq N_0$ (i.e. $h_{Z,V}(N)=P(N)$) and all the higher sheaf cohomologies of $V(N)|_{Z}$ vanish. 
\end{claim}
\begin{proof}
For $[Z]\in \ol{\RHilb_X^{c}(V,P(t))}$ such that $V(N)|_{Z}$ or $V(N)\otimes I_Z$ has higher cohomology, we can increase $N_0$ to kill the higher cohomology, and then apply Noetherian induction. Alternatively, this also follows from \cite[Proposition 4.1]{Dellaca} applied to both $V|_{Z}$ and $V\otimes I_Z$, where $I_Z$ is the ideal sheaf. %There is a globally generated hypothesis required to apply the Proposition, but $V$ is globally generated by assumption and we can apply the usual regularity theorem \cite[Theorem 2.7]{Dellaca} to $I_Z$ and the fact 0-regular implies globally generated \cite[Theorem 1.8.3(i)]{PositivityI}. 
\end{proof}

%Let $H^0(V(N))\cong \mb{A}^M$ as affine spaces. There is a universal family $\mc{S}\rightarrow \mb{A}^M$ such that the fiber over $s\in \mb{A}^M\cong H^0(V(N))$ is $\{s=0\}$. Since $V$ is globally generated, there is a surjection $X\times\mb{A}^M\rightarrow V$ of vector bundles over $X$.
\begin{center}
\begin{tikzcd}
 \mc{Y}\ar[d,"\pi"]\ar[r, "\rho"] & X \\ B:=\left(\ol{\RHilb_X^{c}(V,P(t))}\right)^{\text{red}}& 
\end{tikzcd}
\end{center}
 Let $\pi$ be the universal family restricted to the reduction of the Hilbert scheme and $\rho: \mc{Y}\rightarrow X$ the canonical map that is an embedding on each fiber. By Grauert's theorem, $\pi_{*}\rho^{*}V(N)$ is a vector bundle of rank $P(N)$. Then, we can pull back $H^0(V(N))\otimes \O_X\rightarrow V$ to get $H^0(V(N))\otimes \O_{\mc{Y}}\rightarrow \rho^{*}V$. Pushing forward by $\pi_{*}$ gives us 
\begin{center}
\begin{tikzcd}
H^0(V(N))\otimes \O_{B}\arrow[r] \arrow[bend right=15, rr, "\phi"] & H^0(V(N))\otimes\pi_{*} \O_{\mc{Y}} \arrow[r] & \pi_{*}\rho^{*}V(N)
\end{tikzcd}
\end{center}
%Check: I think the first map is actually an isomorphism by Grauert's theorem
Over a point $[Z]\in B$, $\phi$ is the restriction map $H^0(V(N))\rightarrow H^0(V(N)|_Z)$. From our choice of $N_0$, this restriction map must be surjective. Therefore, $\phi$ is a surjective map of vector bundles and $\ker(\phi)$ is a vector bundle. 

Over a point $[Z]\in B$, $\ker(\phi)$ is precisely the sections of $H^0(V(N))$ that vanish on $Z$. Therefore, if we regard $|\ker(\phi)|$ as the affine bundle corresponding to $\ker(\phi)$, we see $|\ker(\phi)|$ and the reduction of $\wt{D}(V(N),a,P(t))^{\text{cl}}$ agree.
\end{proof}

From Proposition \ref{VBN}, we see 
\begin{corollary}
\label{VBNC}
If $P(t)\in S_c$, 
\begin{align*}
\dim(\wt{D}(V(N),a,P(t))^{\text{cl}}) = p_V(N)-P(N)+\dim(\RHilb_X^{c}(V,a,P(t)))
\end{align*}
for all $N\geq N_0$. Here, $N_0$ depends on $P(t)$. 
\end{corollary}

\begin{corollary}
\label{FYC}
There is $N_0$ dependent on $P(t)$ such that for all $N\geq N_0$, $\ol{D(V(N),a,P(t))}=D(V(N),a,P(t))^{\text{cl}}$. 
\end{corollary}

%Start by implementing referee's comments for this proof...maybe should be broken up into claims?
%Decided: Q the polynomial becomes R
%Q the sheaf stays Q
%replace the dim(X)-rank(V)+a's by c
\Cref{VBN} suffices for proving \Cref{MT}. For the purposes of proving \Cref{incidence} and for the statement in terms of Grothendieck rings (\Cref{GR}), we need \Cref{SC,SI} below. %START HERE

\begin{proposition}
\label{SC}
Given $P(t)\in S_c$, there exists $N_0$ depending on $P(t)$ such that for all $N\geq N_0$, there is a closed subset $E_N\subset \wt{D}(V(N),a,P(t))$, whose codimension in $\wt{D}(V(N),a,P(t))$ is bounded from below by $R(N)$, for $R(t)$ a polynomial with the same degree and leading coefficient as $P(t)$, such that each fiber of
\begin{align*}
\pi_1: \wt{D}(V(N),a,P(t))\backslash E_N\to D(V(N),a,P(t))\backslash \pi_1(E_N)
\end{align*}
is a single reduced point. 
\end{proposition}

In words, the fibers of the map $\wt{D}(V(N),a,P(t))\to D(V(N),a,P(t))$ are single reduced points away from a set of high codimension for large $N$. 

\begin{proof}
%The fibers of $\pi_1: \wt{D}(V(N),a,P(t))\to D(V(N),a,P(t))$ are subschemes of $\RHilb_X^{c}(V,a,P(t))$ by definition, whose dimension is independent of $N$. Therefore, if $A_N\subset D(V(N),a,P(t))$ is a closed subset of codimension $R_1(N)$, then $E(N)=\pi^{-1}(A_N)$ is of codimension $R(N)$ for $R(t)\sim R_1(t)$. 
We will find a suitable $A_N\subset D(V(N),a,P(t))$ and set $E(N)=\pi^{-1}(A_N)$. We will use an incidence correspondence to search for length 2 subschemes in the fibers of $\pi_1: \wt{D}(V(N),a,P(t))\to D(V(N),a,P(t))$ and argue by dimension reasons that fibers containing length 2 subschemes occur in high codimension.

To do this, let $\mc{H}^{[2]}$ be the Hilbert scheme of length two subschemes of $\RHilb_X^{c}(V,a,P(t))$. If $\mc{Z}\to \mc{H}^{[2]}$ is the universal family, there is an induced map $\rho: \mc{Z}\rightarrow \RHilb_X^{c}(V,a,P(t))$. Let $\mc{Y}\to \RHilb_X^{c}(V,a,P(t))$ be the universal family and consider the composite $\rho^{*}\mc{Y}\to \mc{Z}\to \mc{H}^{[2]}$. 
\begin{center}
\begin{tikzcd}
X & \rho^{*}\mc{Y} \arrow[l,swap, "\widetilde{\rho}"] \arrow[d] \arrow[dd, bend right=20, swap, "\phi"] \arrow[r]  & \mc{Y} \arrow[d] \\
   & \mc{Z} \arrow[r, "\rho"] \arrow[d] & \RHilb_X^{c}(V,a,P(t))\\
& \mc{H}^{[2]}
\end{tikzcd}
\end{center}
Now, let $\wt{D}_2(V(N),a,P(t))\subset H^0(V(N))\times \mc{H}^{[2]}$ consist of pairs $(s,[Z])$ such that $s$ pulled back to $(\wt{\rho}^{*}V)|_{\phi^{-1}([Z])}$ vanishes. Equivalently, $s$ vanishes on the scheme-theoretic image of $\phi^{-1}([Z])\to X$. 
\begin{center}
\begin{tikzcd}
& \wt{D}_2(V(N),a,P(t)) \arrow[dl, swap, "\rho_1"] \arrow[dr,"\rho_2"] &\\
D(V(N),a,P(t)) & & \mc{H}^{[2]}
\end{tikzcd}
\end{center}
We want to bound the dimension of a fiber $\rho_2^{-1}([Z])$ for $[Z]\in \mc{H}^{[2]}$. Let $W\subset X$ be the scheme theoretic image of $\phi^{-1}([Z])$ in $X$. We claim $\deg(W)=2\dim(Y)$ for $[Y]\in \RHilb_X^{c}(V,a,P(t))$. If $[Z]$ corresponds to two reduced points in $\RHilb_X^{c}(V,a,P(t))$, then this is clear. If $[Z]$ corresponds to a tangent vector, we apply Lemma \ref{WG}. 

%Notice that we are using the fact that we have restricted $\RHilb_X^{\dim(X)-\rank(V)+a}(V,a,P(t))$ to integral subschemes, otherwise we could for example have a deformation that just moves an embedded point around. Similarly, in the case $[Z]$ corresponds to two reduced points, we could have the varieties above those two points coincide except for embedded points. 

Now, we apply Corollary \ref{TT2}, to see the dimension $\rho_2^{-1}([Z])$ is bounded above by $p_V(N)-R_1(N)$ for a polynomial $R_1(t)$ with the same degree as $P(t)$ and with twice the leading coefficient. Therefore,
\begin{align*}
\dim(\rho_1(\wt{D}_2(V(N),a,P(t)) ))&\leq \dim(\mc{H}^{[2]})+P_V(N)-R_1(N)\\
\dim(\pi_1^{-1}(\rho_1(\wt{D}_2(V(N),a,P(t)) )))&\leq \dim(\mc{H}^{[2]})+P_V(N)-R_1(N)+\dim(\RHilb_X^{c}(V,a,P(t)))
\end{align*}
for all $N>0$. Applying Proposition \ref{VBN} shows that this is less than $\dim(\wt{D}(V(N),a,P(t)))$ for $N$ large, so 
\begin{align*}
\pi_1: \wt{D}(V(N),a,P(t))\to D(V(N),a,P(t))
\end{align*}
is generically injective. Finally, we apply Proposition \ref{VBN} again to see
\begin{align*}
\dim(D(V(N),a,P(t)))-\dim(\pi_1^{-1}(\rho_1(\wt{D}_2(V(N),a,P(t)) )))\geq & \dim(\RHilb_X^{c}(V,a,P(t))) \\&-\dim(\mc{H}^{[2]})-P(N)+R_1(N),
\end{align*}
and the right side is a polynomial with same degree and leading coefficient as $P(t)$, as desired. To finish, we let $E_N$ be the closure of $\rho_1(\wt{D}_2(V(N),a,P(t)))$. 
\end{proof}

Similarly, we also want to bound away the dimension of the intersection $D(V(N),a,P(t))\cap D(V(N),a,R(t))$ for two different polynomials $P(t), R(t)\in S_c$. 
\begin{proposition}
\label{SI}
Given $P(t),R(t)\in S_c$ distinct polynomials, there exists $N_0$ depending on $P(t)$ and $R(t)$ such that for all $N\geq N_0$, the codimension of 
\begin{align*}
D(V(N),a,P(t))\cap D(V(N),a,R(t))\subset H^0(V(N))
\end{align*}
is bounded below by $R_1(N)$ for some polynomial $R_1$ of degree $c$ with leading coefficient the sum of the leading coefficients of $P(t)$ and $R(t)$. 
\end{proposition}
Proposition \ref{SI} is proven in a similar way to Proposition \ref{SC} and without the complication of having to consider tangent vectors in the fibers, so we omit the proof. In this case, we would want to consider the incidence correspondence of triples 
\begin{align*}
(s,[Z_1],[Z_2])\in H^0(V(N))\times \RHilb_X^{c}(V,a,P(t))\times \RHilb_X^{c}(V,a,R(t))
\end{align*}
where $s$ vanishes on both $Z_1$ and $Z_2$. Then, as in the proof of Proposition \ref{SC}, we use the fact that $\deg(Z_1\cup Z_2)=\deg(Z_1)+\deg(Z_2)$ and apply Corollary \ref{TT2}. Equivalently, we could have repeated the argument of Proposition \ref{SC} on 
\begin{align*}
\RHilb_X^{c}(V,a,P(t))\cup \RHilb_X^{c}(V,a,R(t))
\end{align*}
and proved both Proposition \ref{SC} and Proposition \ref{SI} at the same time. 

\subsection{Case of high degree}
\label{HID}
We now deal with the crux of the argument. First, since vector bundles in general don't have a filtration by bundles of smaller rank, we will be naturally led to deal with coherent sheaves that are vector bundles away from a set of high codimension.

\begin{definition}
\label{DAF}
If $F$ is a coherent sheaf on $X$ that is a vector bundle on $X\backslash Z$ for some closed subset $Z\subset X$ of codimension at least $\rank(F|_{X\backslash Z})-a+1$, then $D(F,a)\subset H^0(F)$ is defined to be the locus 
\begin{align*}
\{s\in H^0(F)| \dim(\{s|_{X\backslash Z}=0\})\geq \dim(X)-\rank(F|_{X\backslash Z})+a\}.
\end{align*} 
\end{definition}

Lemma \ref{FOH} shows $D(F,a)$ is closed. 
\begin{remark}
Regarding notation, if $X$ is integral, then instead of saying $F$ is a vector bundle away from a closed subset of small dimension, we could instead define
\begin{align*}
S(F):=\{x\in X| \ \rank(F|_x)>\rank(F)\}
\end{align*}
to be the closed subset where $F$ jumps rank and say $S(F)$ has small dimension. Here, $S(F)$ is the locus where $F$ is not locally free and is called the \emph{singularity set} \cite[Chapter 2, Section 1]{Okonek}. However, since our argument doesn't depend on whether $X$ is reduced, we will instead always refer to a large open set on which $F$ is a vector bundle, instead of just taking the complement of $S(F)$. 
\end{remark}

\subsubsection{A short exact sequence}
\begin{proposition}
\label{SSS}
Let $F$ be a globally generated sheaf on $X$ such that $F|_{X\backslash Z}$ is a vector bundle on $X\backslash Z$ for some proper closed subset $Z\subset X$. If $\codim(Z)\geq \rank(F|_{X\backslash Z})\geq 1$, then there exists a short exact sequence of sheaves
\begin{align*}
0\to \O_X\to F\to Q\to 0,
\end{align*}
where $Q$ is a vector bundle of rank $\rank(F|_{X\backslash Z})-1$ away from some closed subset of codimension at least $\rank(F|_{X\backslash Z})$. 
\end{proposition}

\begin{proof}
The ideas are all in \cite[Example 12.1.11]{Fulton}, but we will describe the minor modifications necessary. Since $F$ is globally generated, we can find a surjection $\O_X^M\to F$ for $M>>0$. If we restrict to $X\backslash Z$, then $\O_{X\backslash Z}^M\to F|_{X\backslash Z}$  is a surjection of vector bundles. Viewing the total spaces of the vector bundles as affine bundles over $X\backslash Z$, we have a surjection of affine bundles $\rho: (X\backslash Z)\times \mb{A}^M\to |(F|_{X\backslash Z})|$, where $|(F|_{X\backslash Z})|$ is the affine bundle over $X\backslash Z$ corresponding to $F|_{X\backslash Z}$. 

Let $W\subset |(F|_{X\backslash Z})|$ be the zero section. Following the argument in \cite[Example 12.1.11]{Fulton}, we can find $t_1,\ldots,t_M\in K$ such that $(T_1-t_1,\ldots,T_M-t_M)$ form a regular sequence on the pullback $\rho^{-1}(W)$. If the surjection $\O_X^M\to F$ is given by sections $s_1,\ldots, s_M\in H^0(F)$, then we choose our map $\O_X\to F$ to be the section $s=t_1s_1+\cdots+t_M s_M$. 

Therefore, $\{s=0\}$ is codimension at least $\rank(F)$. Without loss of generality, we can replace $Z$ with $Z\cup \{s=0\}$ so that $s$ does not vanish on $Z$ and $Z$ is codimension at least $\rank(F)$. Then, the quotient $Q$ is a vector bundle on $X\backslash Z$. 
%This is because locally on open $U\subset X\backslash Z$ small enough, the map $\O_{X\backslash Z}\xrightarrow{s} F_{X\backslash Z}\cong \O_{X\backslash Z}^{\rank(F)}$ is described by $1\rightarrow (v_1,\ldots,v_{\rank(F)})$. At each point of $U$, $v_i$ is not zero for some $i$. We can shrink $U$ and assume $v_1$ is nonwhere zero on $U$. Then, $v_1$ is invertible, and we can change basis, so the map $\O_{X\backslash Z}\xrightarrow{s} F_{X\backslash Z}\cong \O_{X\backslash Z}^{\rank(F)}$ is given by $1\rightarrow (1,0,\ldots,0)$, so $Q$ is a vector bundle on $X\backslash Z$. 
\end{proof}

\subsubsection{Projecting onto the quotient}
The purpose of Section \ref{HID} is to prove the following two propositions, whose proofs are almost identical. We could have also defined the appropriate generalization to Definitions \ref{DI} and \ref{DIp} analogous to how Definition \ref{DAF} generalizes Definition \ref{DA} and proved both results in one statement, but this seemed to make the exposition less clear. 

\begin{proposition}
\label{HD1}
If $F$ is a globally generated sheaf on $X$ restricting to a vector bundle away from some closed subset $Z\subset X$ of codimension at least $\rank(F|_{X\backslash Z})-a+1$, then there exist constants $N_0, N_1$ such that the codimension of $D(F(N),a)\subset H^0(F(N))$ is at least 
\begin{align*}
\frac{1}{(\dim(X)-\rank(F|_{X\backslash Z})+a)!}(N-N_1)^{\dim(X)-\rank(F|_{X\backslash Z})+a}
\end{align*}
for all $N\geq N_0$.%, where $N_0$ depends only on $F$. 
\end{proposition}

\begin{proof}
For the base case, if $\rank(F|_{X\backslash Z})=0$, then $D(F(N),a)=\emptyset$ by definition, so the codimension is infinity according to our convention. This is the key place where we use $a>0$. Now, we will use induction on $\rank(F|_{X\backslash Z})$. 

Apply Proposition \ref{SSS} to get a short exact sequence 
\begin{align*}
0\to \O_X\to F\to Q\to 0,
\end{align*}
where $F$ and $Q$ are vector bundles on $X\backslash Z'$, where $Z'\supset Z$ is codimension at least $\rank(F|_{X\backslash Z})-a+1$ in $X$. Replace $Z$ by $Z'$. 

Pick $N_1'$ to be large enough so that $\O_X(N_1')$ is 0-regular (so we can apply \Cref{TT1}), so in particular the projection $\operatorname{pr}: H^0(F(N))\to H^0(Q(N))$ is surjective for $N\geq N_1'$. We set $N_0=N_1=N_1'$ for now, and will possibly increase $N_0$ and $N_1$ in the remainder of the proof. We write 
\begin{align*}
D(F(N),a)=\operatorname{pr}^{-1}(D(Q(N),a))\cup (D(F(N),a)\backslash \operatorname{pr}^{-1}(D(Q(N),a))). 
\end{align*}
The codimension of $\operatorname{pr}^{-1}(D(Q(N),a))$ in $H^0(F(N))$ is the codimension of $D(Q(N),a)$ in $H^0(Q(N))$, which is at least
\begin{align*}
\frac{1}{(\dim(X)-\rank(F|_{X\backslash Z})+a+1)!}(N-N_1)^{\dim(X)-\rank(F|_{X\backslash Z})+a+1}
\end{align*}
for $N\geq N_0$ (possibly after increasing $N_0, N_1$) by induction on the rank of $F$. 

Therefore, it suffices to bound the codimension of $D(F(N),a)\backslash \operatorname{pr}^{-1}(D(Q(N),a))$. We will do this by bounding the codimension of each fiber of
\begin{align*}
\operatorname{pr}|_{D(F(N),a)\backslash \operatorname{pr}^{-1}(D(Q(N),a))}: D(F(N),a)\backslash \operatorname{pr}^{-1}(D(Q(N),a))\to H^0(Q(N)). 
\end{align*}
Suppose $s\in D(F(N),a)\backslash \operatorname{pr}^{-1}(D(Q(N),a))$. Let $W\subset X$ be the closure of $\{\operatorname{pr}(s)|_{X\backslash Z}=0\}\subset X\backslash Z$ in $X$ equipped with the reduced subscheme structure. By assumption, $\dim(W)=\dim(X)-\rank(F)+a$. Let $W_1,\ldots,W_\ell$ be the components of $W$ of maximal dimension equipped with the reduced subscheme structure. Consider the short exact sequence
\begin{align*}
0\to \O_X(N)\to F(N)\to Q(N)\to 0.
\end{align*}
If $s_1,s_2\in \operatorname{pr}^{-1}(\operatorname{pr}(s))$, then $s_1-s_2$ is in the subspace $H^0(\O_X(N))\subset H^0(F(N))$. If $s_1,s_2$ both vanish on $W_i$, then so does $s_1-s_2$. By Corollary \ref{TT2} applied to the case where the vector bundle is $\O_X(N_1')$ and the degree is 1, we can bound the Hilbert function $h_{W_i,\O_X(N_1')}(N-N_1')$ from below. We see that the locus $\{s'\in \operatorname{pr}^{-1}(\operatorname{pr}(s))| s'|_{W_i}=0\}\subset \operatorname{pr}^{-1}(\operatorname{pr}(s))$ is codimension at least 
\begin{align*}
\binom{N-N_1'+1+(\dim(X)-\rank(F|_{X\backslash Z})+a)}{\dim(X)-\rank(F|_{X\backslash Z})+a},
\end{align*}
which is at least
\begin{align*}
\frac{1}{(\dim(X)-\rank(F|_{X\backslash Z})+a)!}(N-N_1')^{\dim(X)-\rank(F|_{X\backslash Z})+a}.
\end{align*}
% Note, that we are bounding the Hilbert function $h_{W_i,\O_X(N_1')}(N-N_1')$, not the Hilbert function $h_{Z,\O_X}(N)$. 

Repeating this for each $i$ shows 
\begin{align*}
\operatorname{pr}^{-1}(\operatorname{pr}(s))\cap \left(D(F(N),a)\backslash \operatorname{pr}^{-1}(D(Q(N),a))\right)\subset \operatorname{pr}^{-1}(\operatorname{pr}(s))
\end{align*}
is codimension at least $\frac{1}{(\dim(X)-\rank(F|_{X\backslash Z})+a)!}(N-N_1')^{\dim(X)-\rank(F|_{X\backslash Z})+a}$. Applying Lemma \ref{FY6} allows us conclude that the codimension of
\begin{align*}
D(F(N),a)\backslash \operatorname{pr}^{-1}(D(Q(N),a))\subset H^0(F(N))
\end{align*}
is at least $\frac{1}{(\dim(X)-\rank(F|_{X\backslash Z})+a)!}(N-N_1')^{\dim(X)-\rank(F|_{X\backslash Z})+a}$. By construction, $N_1$ is already at least $N_1'$. Finally, we choose $N_0$ large enough so that
\begin{align*}
\frac{1}{(\dim(X)-\rank(F|_{X\backslash Z})+a+1)!}(N-N_1)^{\dim(X)-\rank(F|_{X\backslash Z})+a+1}
\end{align*}
is at least
\begin{align*}
\frac{1}{(\dim(X)-\rank(F|_{X\backslash Z})+a)!}(N-N_1)^{\dim(X)-\rank(F|_{X\backslash Z})+a}
\end{align*}
for $N\geq N_0$. 
\end{proof}

%\begin{Rem}
%Recall at the end of section \ref{NE}, we mentioned that the proof of Proposition \ref{HD1} must use $a>0$ in a crucial way. The proof of Proposition \ref{HD1} uses descending induction on rank, %and the 
%\end{Rem}
\begin{proposition}
\label{HD2}
Fix a degree $d$. There exist constants $N_0, N_1$ such that for all $P(t)$ of degree $c$ and leading coefficient at least $\frac{d}{c!}\rank(V)$ (e.g. $\RHilb_X^{c}(V,a,P(t))$ parameterizes subschemes of degree at least $d$) and for all $N\geq N_0$, the codimension of $D(V(N),a,P(t))$ in $H^0(V(N))$ is at least 
\begin{align*}
\frac{d}{c!}(N-N_1)^{c}.
\end{align*}
\end{proposition}

%Recall that $V$ is assumed to be 0-regular, so it is globally generated \cite[Theorem 1.8.3]{PositivityI}. Also, note that 
The $N_0, N_1$ in the statement of Proposition \ref{HD2} depend on $d$, but do \emph{not} depend on $P(t)$. 

\begin{proof}
We repeat the proof of Proposition \ref{HD1} with $F$ replaced by $V$. Then, when we take the components $W_1,\ldots,W_\ell$ of $W$ of maximal dimension, we can throw out all the components of degree less than $d$. Therefore, when we apply Corollary \ref{TT2}, we get the codimension of $D(V(N),a,P(t))\backslash \on{pr}^{-1}(D(Q(N),a))\subset H^0(V(N))$ is at least $d\binom{N-N_1'-d+c}{c}$. Increasing $N_0, N_1$ if necessary, we can assume $N_0,N_1\geq N_1'+d-1$, so that this is at least 
\begin{align*}
\frac{d}{c!}(N-N_1)^{c}.
\end{align*}
for $N\geq N_0$.

As before, increasing $N_0, N_1$ if necessary, Proposition \ref{HD1} applied to $Q$ shows the codimension of $D(Q(N),a)$ in $H^0(Q(N))$ is at least
\begin{align*}
\frac{1}{(c+1)!}(N-N_1)^{c+1},
\end{align*}
for all $N\geq N_0$. 

Therefore, increasing $N_0$ further implies the codimension of $D(V(N),a,P(t))\subset H^0(V(N))$ is at least $\frac{d}{c!}(N-N_1)^{c}$ for $N\geq N_0$. 
\end{proof}

\subsection{Conclusion of the argument}
We now put together the pieces to finish. 
\begin{comment}
Theorem \ref{DIM} and Theorem \ref{MT} will be implied by Theorem \ref{MT2}.
\begin{theorem}
\label{MT2}
For each $m\geq 0$ and $a>0$, there exists $N_0$ such that for $N\geq N_0$, a component of 
\begin{align*}
D(V(N,a))\backslash \bigcup_{i=1}^{m}\bigcup_{P(t)\sim P_i(t)}D(V(N),a,P(t))
\end{align*}
of maximal dimension is in $\ol{D(V(N),a,P(t))}$ for $P(t)\sim P_{m+1}(t)$. 
\end{theorem}

Corollary \ref{GR} will follow from Theorem \ref{GR2}. Recall the definition of $\widehat{\mc{M}_{\bf{L}}}$ from Section \ref{GRdef}. 

\begin{theorem}
\label{GR2}
In $\widehat{\mc{M}_{\bf{L}}}$, for $m\geq 0$,
\begin{align*}
\lim_{N\to \infty}\frac{[D(V(N),a)\backslash \bigcup_{i=1}^{m}\bigcup_{P(t)\sim P_i(t)}D(V(N),a,P(t))]}{{\bf L}^{P_V(N)-P_{m+1}(N)}} = \sum_{P(t)\sim P_{m+1}(t)}\frac{[\RHilb_{X}^{\dim(X)-\rank(V)+a}]}{{\bf L}^{P(N)-P_{m+1}(N)}}. 
\end{align*}
\end{theorem}
\end{comment}
We will prove Theorems \ref{MT} and \ref{GR} together. 
\begin{proof}
Let $d_{\text{min}}$ be the degree of the varieties parameterized by $\RHilb_X^{c}(V,P_{c,m+1}(t))$, e.g. $d_{\text{min}}$ is the product of $\frac{c!}{\rank(V)}$ and the leading coefficient of $P_{c,m+1}$. Let $d_{\text{med}}=d_{\text{min}}\rank(V)+1$. 

Let $M$ be the integer where $\RHilb_X^{c}(V,P_{c,M}(t))$ parameterizes varieties of degree at least $d_{\text{med}}$, but $\RHilb_X^{c}(V,P_{c,M-1}(t))$ does not. 
%Choose $N_0$ so that the conclusion of Proposition \ref{VBN} hold for all $P(t)$ where $P(t)\sim P_i(t)$ for $i<M$ and we can assume the conclusion of Proposition \ref{SC} holds for $P(t)\sim P_{m+1}(t)$. Increasing $N_0$ if necessary, we can assume the conclusion of Proposition \ref{SI} holds for all choices of $P(t)\neq Q(t)$, where $P(t)\sim P_i(t)$ for $i\leq m$ and $Q(t)\sim P_{m+1}(t)$. Increasing $N_0$ if necessary, we can also assume the conclusion of Proposition \ref{HD2} holds. The choice of $N_0$ in Proposition \ref{HD2} also does not depend on $d$, but we will apply it in the case $d=d_{\text{med}}$. 
Choose $N_0$ so that 
\begin{enumerate}
\item
the conclusion of Proposition \ref{VBN} hold for all $P(t)$ where $P(t)\sim P_{c,i}(t)$ for $i<M$
\item 
the conclusion of Proposition \ref{SC} holds for all $P(t)\sim P_{c,m+1}(t)$
\item
the conclusion of Proposition \ref{SI} holds for all choices of $P(t), Q(t)$, where $P(t)\sim P_{c,i}(t)$ for $i\leq m+1$ and $Q(t)\sim P_{c,m+1}(t)$
\item
the conclusion of Proposition \ref{HD2} holds in the case $d=d_{\text{med}}$. 
\end{enumerate}
%The choice of $N_0$ in Proposition \ref{HD2} also does not depend on $d$, but in any case we will apply it in the case $d=d_{\text{med}}$. 

For the rest of the argument, fix $P(t)\sim P_{c,m+1}(t)$. From Proposition \ref{SC}, we can increase $N_0$ so that $\tilde{D}(V(N),a,P(t))\to D(V(N),a,P(t))$ is birational onto its image for $N\geq N_0$, and from Proposition \ref{SI} we can increase $N_0$ so that 
\begin{align*}
\dim\left(D(V,a,P(t))\backslash \bigcup_{i=1}^{m}\bigcup_{P'(t)\sim P_{c,i}(t)}D(V(N),a,P'(t))\right)=\dim(D(V,a,P(t)))
\end{align*}
for $N\geq N_0$. Then, if $N\geq N_0$ and $Q(t)\sim P_{c,i}(t)$ for $m+1<i<M$,
\begin{align*}
\dim\left(D(V,a,P(t))\backslash \bigcup_{i=1}^{m}\bigcup_{P'(t)\sim P_{c,i}(t)}D(V(N),a,P'(t))\right)-\dim(D(V,a,Q(t)))
\end{align*}
is bounded below by a nonconstant polynomial with positive leading coefficient in $N$.  Similarly, Proposition \ref{HD2} shows that the same statement is true for $Q(t)\sim P_{c,i}(t)$ for $i\geq M$. This shows Theorem \ref{MT}.

To see Theorem \ref{GR}, we note in addition by Proposition \ref{SI} that 
\begin{align*}
\dim(D(V,a,P(t)))-\dim(D(V,a,P(t))\cap D(V(N),a,Q(t)))
\end{align*}
is bounded below by a nonconstant polynomial with positive leading coefficient for $N\geq N_0$ if $P(t)\neq Q(t)\sim P_{c,i}(t)$ for $i\leq m+1$, and by Proposition \ref{SC} there is a set $E(N)\subset \wt{D}(V,a,P(t))$ such that the fibers of the map
\begin{align*}
\wt{D}(V,a,P(t))\backslash E(N)\to D(V,a,P(t))
\end{align*}
are either empty or a single reduced point, and the codimension of $E(N)$ in $\wt{D}(V,a,P(t))$ is bounded below by a nonconstant polynomial with positive leading coefficient. To conclude we apply Proposition \ref{ISOG}. 
\end{proof}

\appendix
\section{Facts about Hilbert functions}
Throughout the appendix, we assume that we are working in a projective scheme $X$ and $V$ is a vector bundle on $X$. 
\subsection{Lower bound on Hilbert function given degrees}
To give crude bounds on Hilbert function, we will generalize a well-known lemma (see \cite[Lemma 3.1]{Montreal}). Recall $h_{Z,E}(n):= \dim(\im(H^0(V(n))\rightarrow H^0(V(n)|_{Z})))$ (\Cref{HIF}).
\begin{lemma}
\label{M1}
If $H^1(V(n-1))=0$, $Z\subset X$ is a closed subscheme and $H$ is a hyperplane section of $X$ that is a nonzero divisor when restricted to $Z$, then
\begin{align*}
h_{Z,V}(n)\geq h_{Z,V}(n-1)+h_{Z\cap H,V}(n).
\end{align*}
\end{lemma}

\begin{proof}
Consider the diagram
\begin{center}
\begin{tikzcd}
0 \arrow[r] & H^0(V(n-1)) \arrow[r, "\times H"] \arrow[d]& H^0(V(n)) \arrow[r] \arrow[d]& H^0(V(n)|_H)  \arrow[r] \arrow[d]& 0\\
0\arrow[r]  & H^0(V(n-1)|_{Z}) \arrow[r, "\times H"] & H^0(V(n)|_{Z}) \arrow[r] & H^0(V(n)|_{H\cap Z}) &
\end{tikzcd}
\end{center} 
The rows are exact on the left since $H$ is a nonzero divisor on $X$ and $Z$. Taking the image under the middle vertical map yields a short exact sequence with $\im(H^0(V(n))\to H^0(V(n)|_{Z}))$ in the middle. 
    \begin{align*}
    0 \to \im(\im(H^0(V(n))\rightarrow H^0(V(n)|_{Z}))\cap \im(H^0(V(n-1)|_{Z})\xrightarrow{\times H} H^0(V(n)|_{Z})))\to\\ \im(H^0(V(n))\to H^0(V(n)|_{Z})) \to \im(H^0(V(n))\to H^0(V(n)|_{H\cap Z}))\to  0
    \end{align*}
Then,
\begin{align*}
h_{Z,V}(n):=& \dim(\im(H^0(V(n))\rightarrow H^0(V(n)|_{Z})))\\
=&\dim(\im(H^0(V(n))\rightarrow H^0(V(n)|_{Z}))\cap \im(H^0(V(n-1)|_{Z})\rightarrow H^0(V(n)|_{Z})))+\\
& \dim(\im(H^0(V(n))\rightarrow H^0(V(n)|_{H\cap Z})))\\
\geq & \dim(\im(H^0(V(n-1))\rightarrow H^0(V(n-1)|_{Z})))+\dim(\im(H^0(V(n))\rightarrow H^0(V(n)|_{H\cap Z})))\\
=& h_{Z,V}(n-1)+h_{Z\cap H,V}(n).
\end{align*}
\end{proof}

\begin{comment}
\begin{proof}
This is Lemma 3.1 in \cite{Montreal}. Since our field is infinite, such an $H$ exists. Let $E_n\subset H^{0}(X,\mathscr{O}_X(n))$ be the image of $H^{0}(\mb{P}^r,\mathscr{O}_{\mb{P}^r}(n))$. We want to bound $\dim(E_n)-\dim(E_{n-1})$ by $h_{X\cap H}(n)$.

We have the map $E_{n}\rightarrow H^{0}(X\cap H, \mathscr{O}_{X\cap H}(n))$ and we can let $E_{n}'$ be the kernel. We see that $\dim(E_n)-\dim(E_n')=h_{X\cap H}(n)$. Also, we see that, we have a map $E_{n-1}\hookrightarrow E_n$ given by multiplying by $H$, and the image $H\cdot E_{n-1}$ of the map is contained in $E_{n}'$. Therefore, $\dim(E_n)-\dim(E_{n-1})\geq h_{X\cap H}(n)$, and equality holds if and only if $H\cdot E_{n-1}$ is equal $E_{n}'$.
\end{proof}
\end{comment}

\begin{definition}
A vector bundle $E$ on $X$ is said to be $k$-very ample for $k$ a nonnegative integer if the restriction map $H^0(E)\to H^0(E|_{Z})$ is surjective for every 0-dimensional length $k+1$ subscheme $Z\subset X$. 
\end{definition}

The following lemma should be well-known to experts, but the author couldn't find it in the literature.
\begin{lemma}
\label{GGKVA}
If $V$ is globally generated and $L$ is a $k$-very ample line bundle, then $V \otimes L$ is $k$-very ample.
\end{lemma}

\begin{proof}
Let $Z\subset X$ be a length $k+1$ subscheme and $Z=\Spec(A)$ for an Artinian ring $A$. Let $A=\oplus_{i=1}^{\ell}{A_\ell}$, where each $A_i$ is a local Artinian ring. Consider
\begin{center}
\begin{tikzcd}
H^0(V\otimes L) \arrow[r, "\phi"] \arrow[rr, bend left=20, "\psi"] & H^0((V\otimes L )|_{Z}) \arrow[r,"\sim"] & A^{\rank(V)} \\
H^0(V)\otimes H^0(L) \arrow[u] \arrow[r] & H^0(V|_{Z})\otimes H^0(L|_{Z}) \arrow[u] \arrow[r,"\sim"] &  A^{\rank(V)}\otimes_A A \arrow[u] & 
\end{tikzcd}
\end{center}
To show $\phi$ is surjective, it suffices to show that $\psi$ hits every element $e_{ij}\in A^{\rank(V)}\cong \oplus_{i=1}^{\ell}{A_{\ell}}$, where $e_{ij}$ is the element that projects to 0 under each projection $ A^{\rank(V)}\to A_{p}^{\rank(V)}$ for $p\neq i$ and maps to the element  
\begin{align*}
(0,0,\ldots, \underbrace{1}_{\text{index}\ j}, \ldots, 0)\in A_i^{\rank(V)}
\end{align*}
for $p=i$. For convenience of notation, assume $j=1$. Let $m_i\subset A_i$ be the maximal ideal. Since $V$ is globally generated, we can find $s\in H^0(V)$ such that $s$ maps to $(1,0,\ldots,0)$ under the composition
\begin{align*}
H^0(V)\to H^0(V|_Z)\xrightarrow{\sim} A^{\rank(V)}\to A_i^{\rank(V)}\to (A_i/m_i)^{\rank(V)}.
\end{align*}
Let $a\in A$ be the image of $s$ in $A_i$ under the composition 
\begin{align*}
H^0(V)\to H^0(V|_Z)\xrightarrow{\sim} A^{\rank(V)}\to A_i^{\rank(V)}\to A_i,
\end{align*}
where the last map is projection onto the first factor. In particular, $a$ is invertible in $A_i$. Since $L$ is $k$-very ample, the composition
\begin{align*}
H^0(L)\to H^0(L|_Z)\xrightarrow{\sim} A
\end{align*}
is surjective, so we can find $t\in H^0(L)$ such that $t$ projects to 0 in each factor $A_p$ for $p\neq i$ and $t$ projects to $a^{-1}$ in the factor $A_i$. Then, the element $s\otimes t\in H^0(V)\otimes H^0(L)$ maps to a section of $H^0(V\otimes L)$ that maps to $e_{ij}$ under $\psi$ above. 
\end{proof}

\begin{lemma}
\label{TT1}
If $V$ 0-regular in the sense of Castelnuovo-Mumford regularity, $Z\subset X$ a subscheme of degree at least $d$ and dimension $m$, then for $n\geq 1$
\begin{align*}
h_{Z,V}(n-1)\geq \begin{cases}
\rank(V)\binom{n+m}{m+1} &\text{ if $n\leq d$}\\
\rank(V)\sum_{i=0}^{m}{\binom{n-d+i-1}{i}\binom{d+m-i}{m-i+1}} &\text{ if $n> d$}.
\end{cases}
\end{align*}
\end{lemma}

\begin{proof}
We will use induction and Lemma \ref{M1}. First, we show the base case $m=0$. Let $Z=\Spec(A)$ for an Artinian $K$-algebra $A$ and fix an trivialization $V|_{Z}\cong \O_{\Spec(A)}^{\rank(V)}$. We want to show $h_{Z,V}(n-1)\geq \rank(V)\min\{n,d\}$. Since $V$ is globally generated \cite[Theorem 1.8.3]{PositivityI} and $\O(n-1)$ is $n$-very ample \cite[Theorem 1.1]{kVeryAmple}, we can apply Lemma \ref{GGKVA} to see $V(n-1)$ is $n$-very ample. 

Without loss of generality, we can assume $d\leq n$, or else we can just restrict $Z$ to a subscheme of length $n$. (Recall $A$ has quotients of any length less than $A$. To show this, it suffices to see that if $A$ is local, then there is an ideal of $A$ of length 1. To exhibit such an ideal, we can pick any element $a\in A$ that is anniliated by the maximal ideal of $A$ and take the ideal generated by $a$.) Then, by $n$-very ampleness, we see $H^0(V(n-1))\to H^0(V(n-1)|_{Z})$ is surjective, so $h_{Z,V}(n-1)=n$. The reason why we only prove an inequality is that we restricted $Z$ above. 

Next, for the induction step, suppose we know Lemma \ref{TT1} in dimension $m-1$. From the long exact sequence, we have the exactness of
\begin{align*}
H^i(V(-i))\to H^i(V(-i)|_{H})\to H^{i+1}(V(-i-1))
\end{align*}
so $V|_{H}$ is still 0-regular. Also, $H^0(V)\rightarrow H^0(V|_{H})$ is surjective as $H^1(V(-1))=0$. Therefore, we can apply the induction hypothesis and Lemma \ref{M1}.

If $Z$ is dimension $m$ and degree $d$, we have for $n\leq d$
\begin{align*}
h_{Z,V}(n-1)&\geq  h_{Z,V}(-1)+h_{Z\cap H,V}(0)+h_{Z\cap H,V}(1)+\cdots+ h_{Z\cap H,V}(n-1)\\
&\geq 0+\rank(V)\binom{m}{m}+\cdots+\rank(V)\binom{m+n-1}{m}=\rank(V)\binom{n+m}{m+1}.
\end{align*}
If $n>d$,
\begin{align*}
h_{Z,V}(n-1)&\geq h_{Z,V}(d-1)+\left( h_{Z\cap H,V}(d)+\cdots+h_{Z\cap H,V}(n-1)\right)\\
&= \rank(V)\binom{d+m}{m+1}+\sum_{j=d+1}^{n}{\rank(V)\sum_{i=0}^{m-1}{\binom{j-d+i-1}{i}\binom{d+m-1-i}{m-1-i+1}}}\\
&= \rank(V)\binom{d+m}{m+1}+\rank(V)\sum_{i=0}^{m-1}{\binom{d+m-1-i}{m-1-i+1}\sum_{j=d+1}^{n}{\binom{j-d+i-1}{i}}}\\
&= \rank(V)\binom{d+m}{m+1}+\rank(V)\sum_{i=0}^{m-1}{\binom{d+m-1-i}{m-1-i+1}\binom{n-d+i}{i+1}}\\
&= \rank(V)\binom{d+m}{m+1}+\rank(V)\sum_{i=1}^{m}{\binom{d+m-i}{m-i+1}\binom{n-d+i-1}{i}}.
\end{align*}
\end{proof}

\begin{corollary}
\label{TT2}
If $Z\subset X$ is a subscheme of degree at least $d$ and dimension $m$, then for $n\geq 0$,
\begin{align*}
h_{Z,V}(n)\geq d\cdot \rank(V)\binom{n-d+m}{m}
\end{align*}
for a $0$-regular vector bundle $V$ on $X$. 
\end{corollary}
\begin{proof}
We bound $h_{Z,V}(n)$ from below by the $i=m$ term in summation in the statement of Lemma \ref{TT1}. 
\end{proof}

\section{Scheme theoretic facts}
\begin{lemma}
\label{FOH}
Suppose $X$ is an integral projective scheme and $Z\subset X$ a closed subset. Let $S$ be a finite type $K$-scheme, $Y\subset S\times (X\backslash Z)$ be a closed subset and $\pi:Y\rightarrow S$ be the projection. Then,
\begin{align*}
\{s\in S: \dim(\pi^{-1}(s))\geq d\}
\end{align*}
is a closed subset of $S$ for all $d>\dim(Z)$.
\end{lemma}

\begin{proof}
Let $\overline{Y}$ be the closure of $Y$ in $S\times X$ and $\pi:\overline{Y}\rightarrow S$ be the projection.
Let $Y_d\subset Y$ be the closed subset $\{y\in Y: \dim(\pi^{-1}(\pi(y))\cap Y)\geq d\}$. Let $\overline{Y}_d\subset \overline{Y}$ be the closed subset $\{y\in \overline{Y}: \dim(\pi^{-1}(\pi(y)))\geq d\}$.

Since $Y_d\subset \overline{Y}_d$, $\pi(Y_d)\subset \pi(\overline{Y}_d)$. We want to show $\pi(Y_d)= \pi(\overline{Y}_d)$. 

To see the other inclusion, we first claim $Y_d=\overline{Y}_d\backslash (S\times Z)$. Given this, $\pi(Y_d)\supset \pi(\overline{Y}_d)$ also holds, as if $s\in \pi(\overline{Y}_d)$ is a closed point, then $\pi^{-1}(s)$ has dimension at least $d$ in $X$, so $\pi^{-1}(s)\backslash Z$ still has dimension at least $d$ as $\dim(Z)<d$.

To see $Y_d=\overline{Y}_d\backslash (S\times Z)$, suppose $\xi\in \overline{Y}$ is a scheme theoretic point with $\dim(\overline{\{\xi\}})\geq d$ that maps to a closed point of $S$ and $y\in \overline{\{\xi\}}\backslash Z$. Equivalently, $y\in \overline{Y}_d\backslash (S\times Z)$. But then $\overline{\{\xi\}}$ cannot be contained in $Z$ since $\dim(Z)<\dim(\overline{\{\xi\}})$, so $\xi\in Y$ and $y\in Y_d$.
\end{proof}

\begin{comment}
\begin{Lem}
[{\cite[Lemma 7.1]{ACGH2}}]
\label{CFF}
Given a scheme $A$ and two closed subschemes $B,C\subset \mb{P}^r_A$ with $B$ flat over $A$, there exists a closed subscheme $D\subset A$ such that any morphism $T\rightarrow A$ factors through $D$ if and only if $B\times_A T$ is a subscheme of $C\times_A T$.
\end{Lem}
\end{comment}

\begin{lemma}
\label{WG}
Suppose $Z\subset \mb{P}^r$ is a integral projective scheme, and $\mc{Z}\subset \mb{P}^r_{\Spec(K[\epsilon]/(\epsilon^2))}$ is a nontrivial infinitesimal deformation of the embedding $Z\subset \mb{P}^r$. Let $Z'\subset \mb{P}^r$ be the scheme theoretic image of $\mc{Z}\rightarrow \mb{P}^r$. Then, $\deg(Z')=2\deg(Z)$. 
\end{lemma}

\begin{proof}
Assume $Z$ is not contained in hyperplane $\{X_0=0\}$. In the chart $\{X_0\}$, the data of the embedding $Z\subset \mb{P}^r$ is the surjection $K[x_1,\ldots,x_r]\rightarrow A$ for $\Spec(A)=Z\backslash \{X_0=0\}$ and similarly for $\mc{Z}\subset \mb{P}^r_{\Spec(K[\epsilon]/(\epsilon^2))}\rightarrow \mc{A}$ for $\Spec(\mc{A})=\mc{Z}\backslash\{X_0=0\}$. 
\begin{center}
\begin{tikzcd}
K[x_1,\ldots,x_r] \arrow[r, shift left] \arrow[two heads, d] & K[\epsilon,x_1,\ldots,x_r]/(\epsilon^2) \arrow[l,shift left] \arrow[d, two heads]\\
A & \mc{A} \arrow[l] 
\end{tikzcd}
\end{center}
The scheme theoretic image of $\mc{Z}\rightarrow \mb{P}^r$ restricted to $\{X_0=0\}$ can be computed affine locally \cite[Tag 01R8]{stacks-project} and is $\Spec(B)$ for $B=\im(K[x_1,\ldots,x_r]\rightarrow \mc{A})$. Equivalently, it is defined by the ideal $I\subset K[x_1,\ldots,x_r]$, where $I=\ker(K[x_1,\ldots,x_r]\rightarrow \mc{A})$. Let $\eta\in \Spec(A)$ be the generic point. It suffices to show the multiplicity of $B$ at $\eta$ is 2, where multiplicity is defined to be the length of $B_{\eta}$ \cite[Section 1.5]{Fulton}. In the argument below, we will use $A$ is integral to see $A\rightarrow A_{\eta}$ is injective. 
\begin{center}
\begin{tikzcd}
0 \arrow[r] & I \arrow[r] \arrow[d] & K[x_1,\ldots,x_r] \arrow[r] \arrow[d] & B \arrow[r] \arrow[d] & 0\\
0 \arrow[r] & I_{\eta} \arrow[r] & K[x_1,\ldots,x_r]_{\eta} \arrow[r] & B_{\eta} \arrow[r] & 0
\end{tikzcd}
\end{center}
From exactness of localization, it suffices to show that the induced map $B_{\eta}\rightarrow A_{\eta}$ is not an isomorphism. We will assume $B_{\eta}\rightarrow A_{\eta}$ is an isomorphism, and we will show that the deformation $\Spec(\mc{A})\rightarrow \Spec(K[\epsilon]/(\epsilon^2))$ must be trivial. 
\begin{center}
\begin{tikzcd}
0 \arrow[r] & \epsilon A_{\eta} \arrow[r] & \mc{A}_{\eta} \arrow[r] & A_{\eta} \arrow[r] & 0\\
0 \arrow[r] & \epsilon A \arrow[r] \arrow[u, hook] & \mc{A}\arrow [r] \arrow[u] & A \arrow[r] \arrow[u, hook] & 0\\
0 \arrow[r] & \epsilon K \arrow[r] \arrow[u] & K[\epsilon]/(\epsilon^2) \arrow[r] \arrow[u] & K\arrow[r] \arrow[u] & 0
\end{tikzcd}
\end{center}
The middle row is exact by flatness \cite[Proposition 2.2]{Deformation}. From the diagram, we see $\mc{A}\rightarrow \mc{A}_{\eta}$ is an injection. 
\begin{center}
\begin{tikzcd}
B \arrow[r, hook] \arrow[rr,  two heads, bend left] \arrow[d] & \mc{A} \arrow[d, hook] \arrow[r] & A \arrow[d, hook] \\
B_{\eta} \arrow[r] \arrow[rr, "\sim", bend right] & \mc{A}_{\eta} \arrow[r] & A_{\eta}
\end{tikzcd}
\begin{tikzcd}
K[\epsilon,x_1,\ldots,x_r]/(\epsilon^2) \arrow[r, two heads] \arrow[d, shift left] & \mc{A} \arrow[d, shift left] & \\
K[x_1,\ldots,x_r] \arrow[r, two heads] \arrow[u, shift left] & A \arrow[u, shift left] \arrow[r, shift left] & A\otimes_{K} K[\epsilon]/(\epsilon^2) \arrow[dotted, ul] \arrow[l, shift left]
\end{tikzcd}
\end{center}
Since $B\to \mc{A}_{\eta}$ is injective, we know $B\to B_\eta$ is injective. Then, $B\to A_\eta$ is injective, so $B\to A$ is injective. But $B\rightarrow A$ is surjective, so $B\rightarrow A$ is an isomorphism. This gives a splitting $A\rightarrow \mc{A}$ compatible with the splitting $K[x_1,\ldots,x_r]\rightarrow K[\epsilon, x_1,\ldots,x_r]/(\epsilon^2)$. This shows that $\Spec(\mc{A})\rightarrow \Spec(K[\epsilon]/(\epsilon^2))$ is the trivial deformation. 
\end{proof}

\begin{comment}
\begin{Lem}
\label{RF}
If $W\rightarrow Z$ is a finite map of finite type $K$-schemes where each fiber is either empty or a single reduced point, then $W\rightarrow Z$ is an immersion. If particular, if $W\to Z$ is surjective and $W$ is reduced, then $W\to Z$ is an isomorphism. 
\end{Lem}

\begin{proof}
Since finite implies affine, we can reduce to the case where we have a map $A\rightarrow B$ of rings. Our assumption on fibers means, for each maximal ideal $m$ of $A$, the induced map $A/m\rightarrow B/m$ is surjective. Since $B$ is finite over $A$, the cokernel $M$ is a finitely generated $A$-module. Tensoring with $A\rightarrow B\rightarrow M\rightarrow 0$ with $A/m$ yields $A/m\rightarrow B/mB\rightarrow M/mM\rightarrow 0$, so $M/mM=0$. By Nakayama's lemma, this means $M_{m}=0$. Since this is true for all $M$, $M=0$ and $A\rightarrow B$ is surjective.
\end{proof}
\end{comment}

\begin{comment}
\begin{Lem} [{{\cite[Tag 01R8]{stacks-project}}}]
\label{IHML}
Suppose $f:X\rightarrow Y$ is a quasicompact map of schemes with scheme theoretic image $Z$. Then, for open $U\subset Y$, the scheme theoretic image of $f|_{f^{-1}(U)}$ is $Z\cap U$.
\end{Lem}
\end{comment}

\section{Isomorphism in the Grothendieck ring}

%See \cite[Lemma 2.1]{Beke} for a similar statement to Proposition \ref{ISOG}. In particular, s
We will need a criterion for isomorphism in the Grothendieck ring.
\begin{proposition}
\label{ISOG}
Suppose $f: Y\rightarrow Z$ is a morphism of finite type $K$-schemes and $A\subset Y$ and $B\subset Z$ are constructible sets such that $f(A)=B$ and $f^{-1}(p)$ is a single reduced point for each closed $p\in B$. Then $[A]=[B]$ in the Grothendieck ring $\mc{M}$. 
\end{proposition}

\begin{remark}
In characteristic 0, the conclusion of Proposition \ref{ISOG} holds if $f:A\to B$ is only a bijection of closed points \cite[Lemma 2.1]{Beke}. As a note of caution, we have contacted the author and Lemma 2.1 in \cite{Beke} is correct in characteristic zero, but in characteristic $p$ there is a gap in the proof, since separable maps do not necessarily remain separable when restricted to closed subsets. In our case this is guaranteed by the assumption that all of our fibers are reduced points. 
\end{remark}

\begin{proof}
By assumption $f$ induces a bijection $A\to B$. If we write $B=U_1\cup \cdots U_{\ell}$ as a finite union of locally closed sets, it suffices to show $[f^{-1}(U_i)]=[U_i]$. Therefore, it suffices to show the case where $B$ is locally closed. Similarly, we can break $B$ up further, so that $B$ is irreducible and locally closed. 

Since $B$ is locally closed, $A$ is also locally closed. We can equip $A$ and $B$ with the reduced subscheme structure of an open subscheme. Then, $f: A\to B$ is a map of irreducible finite type $K$-schemes such that each fiber over a closed point $p\in B$ is a single reduced point. Pick an affine open $\Spec(R_B)\subset B$ and an affine open $\Spec(R_A)\subset f^{-1}(\Spec(R_B))$. 

By Grothendieck's generic freeness lemma \cite[Theorem 14.4]{View}, we can restrict $R_B$ so that $R_A$ is free over $R_B$. By the assumption on the fibers of $f: A\to B$, $R_A$ must be rank at most 1 over $R_B$, so the map $R_B\to R_A$ is an isomorphism. Then, $\Spec(R_A)\to \Spec(R_B)$ is an isomorphism and we use Noetherian induction on the complement $f: A\backslash \Spec(R_A)\to B\backslash \Spec(R_B)$. 
\end{proof}

\begin{figure}[ht]
    \centering
\begin{longtable}{l | p{11cm} | r}
Term & Informal Definition & Reference \\\hline
$K$ & our base field, which is of arbitrary characteristic and algebraically closed & \Cref{defsec}\\
$Z^{\operatorname{red}}$ & the reduction of the scheme $Z$ & \Cref{defsec}\\
$\Hilb_X$ & Hilbert scheme of subschemes of $X$ & \Cref{hilbdef}\\
$\Hilb_X^c$ & components of $\Hilb_X$ consisting of dimension $c$ subschemes & \Cref{hilbdef}\\
$\RHilb_X$ & open locus in $\Hilb_X$ consisting of (geometrically) integral subschemes & \Cref{hilbdef}\\
$\RHilb^c_X$ & components of $\RHilb_X$ consisting of dimension $c$ subschemes & \Cref{hilbdef}\\
$\RHilb_X^{c}(V,P(t))$ & components of $\RHilb_X$ consisting of subschemes $Z$ with Hilbert polynomial $\chi(V|_{Z}(t))=P(t)$ & \Cref{hilbdefP}\\
$\ol{\RHilb_X^{c}(V,P(t))}$ & closure of $\RHilb_X^{c}(V,P(t))$ in $\Hilb_X$ & \Cref{hilbdefPc}\\ 
$S_c$ & the polynomials $P(t)$ for which $\RHilb_X^{c}(V,P(t))$ is nonempty & \Cref{Sc}\\
$P_V(t)$ & Hilbert polynomial of $V$ & \Cref{pV}\\
$h_{Z,E}$ & Hilbert function given by $\dim(\im(H^0(E(n))\rightarrow H^0(E(n)|_{Z})))$ & \Cref{HIF}\\
$P$ dominates $Q$&$P(t)-Q(t)$ asymptotically approaches $\infty$ & \Cref{PQ}\\
$P\sim Q$ & $P(t)-Q(t)$ is a constant & \Cref{PQ}\\
$P_{c,1},P_{c,2},\ldots $ & A choice of one representative of all the Hilbert polynomials of $\RHilb_X^c$ under equivalence $\sim$ in increasing order by dominance & \Cref{Pi}\\
$D(V,a)$ & sections of $V$ vanishing to dimension $a$ greater than expected & \Cref{DA}\\
$D(F,a)$ & If $F$ is a vector bundle $V$, then this is $D(V,a)$. Here, we allow $F$ to be a coherent sheaf, but it needs to be a vector bundle away from a set of codimension at least $\rank(F)-a+1$ & \Cref{DAF}\\
$\wt{D}(V, a, P(t))$ & $(s,[Z])\in H^0(V)\times \RHilb_X^{\dim(X)-\rank(V)+a}(V,P(t))$, where $s$ vanishes on $Z$ & \Cref{DI}\\
$D(V,a,P(t))$ & sections in $D(V,a)$ vanishing on a (geometrically) integral subscheme of $X$ of dimension $\dim(X)-\rank(V)+a$ and Hilbert polynomial $P(t)$ & \Cref{DIp}\\
$\wt{D}(V, a, P(t))^{\text{cl}}$ & a compactification of $\wt{D}(V, a, P(t))$ where $Z$ is now allowed to be in the closure of $\RHilb_X^{\dim(X)-\rank(V)+a}(V,P(t))$ & \Cref{DIc}\\
$D(V,a,P(t))^{\text{cl}}$ & a compactification of $D(V,a,P(t))$ where the sections are allowed to vanish on a subscheme in the closure of the (geometrically) integral subschemes & \Cref{DIc}\\
$\ol{D(V,a,P(t))}$ & The closure of $D(V,a,P(t))$ in the affine space $H^0(V)$. This will agree with $D(V,a,P(t))^{\text{cl}}$ after replacing $V$ with a high enough twist (see \Cref{FYC}) &\\
$\mc{M}$ & Grothendieck ring of varieties & \Cref{GRdef}\\
$\mb{L}$ & $[\mb{A}^1]\in\mc{M}$ & \Cref{GRdef}\\
$\mc{M}_{\bf{L}}$ & the commutative ring $\mc{M}$ with $\bf{L}$ inverted & \Cref{GRdef}\\
$\widehat{\mc{M}_{\bf{L}}}$ & the completion of $\mc{M}_{\bf{L}}$ with respect to the filtration of ideals given by dimension & \Cref{GRdef}\\\hline
\end{longtable}
\caption{Table of symbols}
\label{symboltable}
\end{figure}

\bibliographystyle{plain}
\bibliography{vbref}

\begin{thebibliography}{10}

\bibitem{Beke}
Tibor Beke.
\newblock The {G}rothendieck ring of varieties and of the theory of
  algebraically closed fields.
\newblock {\em J. Pure Appl. Algebra}, 221(2):393--400, 2017.

\bibitem{Bourqui}
David Bourqui.
\newblock Asymptotic behaviour of rational curves.
\newblock {\em preprint}, 2011.
\newblock arXiv:1107.3824v1.

\bibitem{CNS11}
Raf Cluckers, Johannes Nicaise, and Julien Sebag, editors.
\newblock {\em Motivic integration and its interactions with model theory and
  non-{A}rchimedean geometry. {V}olume {I}}, volume 383 of {\em London
  Mathematical Society Lecture Note Series}.
\newblock Cambridge University Press, Cambridge, 2011.

\bibitem{Dellaca}
Roger Dellaca.
\newblock Gotzmann regularity for globally generated coherent sheaves.
\newblock {\em J. Pure Appl. Algebra}, 220(4):1576--1587, 2016.

\bibitem{EGA}
Jean Dieudonn{\'e} and Alexander Grothendieck.
\newblock \'{E}l\'ements de g\'eom\'etrie alg\'ebrique.
\newblock {\em Inst. Hautes \'Etudes Sci. Publ. Math.}, 4, 8, 11, 17, 20, 24,
  28, 32, 1961--1967.

\bibitem{View}
David Eisenbud.
\newblock {\em Commutative algebra}, volume 150 of {\em Graduate Texts in
  Mathematics}.
\newblock Springer-Verlag, New York, 1995.
\newblock With a view toward algebraic geometry.

\bibitem{Fulton}
William Fulton.
\newblock {\em Intersection theory}, volume~2.
\newblock Springer-Verlag, Berlin, second edition, 1998.

\bibitem{Montreal}
Joe Harris.
\newblock {\em Curves in projective space}, volume~85 of {\em S\'eminaire de
  Math\'ematiques Sup\'erieures [Seminar on Higher Mathematics]}.
\newblock Presses de l'Universit\'e de Montr\'eal, Montreal, Que., 1982.
\newblock With the collaboration of David Eisenbud.

\bibitem{Deformation}
Robin Hartshorne.
\newblock {\em Deformation theory}, volume 257 of {\em Graduate Texts in
  Mathematics}.
\newblock Springer, New York, 2010.

\bibitem{kVeryAmple}
Yukitoshi Hinohara, Kazuyoshi Takahashi, and Hiroyuki Terakawa.
\newblock On tensor products of {$k$}-very ample line bundles.
\newblock {\em Proc. Amer. Math. Soc.}, 133(3):687--692, 2005.

\bibitem{Kollar}
J\'anos Koll\'ar.
\newblock {\em Rational curves on algebraic varieties}, volume~32.
\newblock Springer-Verlag, Berlin, 1996.

\bibitem{PositivityI}
Robert Lazarsfeld.
\newblock {\em Positivity in algebraic geometry. {I}}, volume~48.
\newblock Springer-Verlag, Berlin, 2004.
\newblock Classical setting: line bundles and linear series.

\bibitem{Liu}
Qing Liu.
\newblock {\em Algebraic geometry and arithmetic curves}, volume~6 of {\em
  Oxford Graduate Texts in Mathematics}.
\newblock Oxford University Press, Oxford, 2002.
\newblock Translated from the French by Reinie Ern\'{e}, Oxford Science
  Publications.

\bibitem{Okonek}
Christian Okonek, Michael Schneider, and Heinz Spindler.
\newblock {\em Vector bundles on complex projective spaces}.
\newblock Modern Birkh\"auser Classics. Birkh\"auser/Springer Basel AG, Basel,
  2011.
\newblock Corrected reprint of the 1988 edition, With an appendix by S. I.
  Gelfand.

\bibitem{stacks-project}
The {Stacks Project Authors}.
\newblock {\itshape Stacks Project}.
\newblock \url{http://stacks.math.columbia.edu}, 2017.

\bibitem{TRockyCIH}
Dennis Tseng.
\newblock A note on $r$ hypersurfaces intersecting in $\mathbb{P}^r$.
\newblock {\em preprint}.
\newblock arXiv:1908.01620.

\bibitem{Tseng}
Dennis Tseng.
\newblock {Collections of Hypersurfaces Containing a Curve}.
\newblock {\em International Mathematics Research Notices}, 06 2018.

\bibitem{VakilWood}
Ravi Vakil and Melanie~Matchett Wood.
\newblock Discriminants in the {G}rothendieck ring.
\newblock {\em Duke Math. J.}, 164(6):1139--1185, 2015.

\end{thebibliography}
\end{document}